\documentclass[journal]{new-aiaa}
\usepackage[utf8]{inputenc}

\usepackage{graphicx}
\usepackage{amsmath}
\usepackage[version=4]{mhchem}
\usepackage{siunitx}
\usepackage{longtable,tabularx}
\usepackage{algorithm}
\usepackage{algorithmic}
\title{Deep Networks as Approximators of Optimal Transfers Solutions in Multitarget Missions}

\author{Haiyang Li \footnote{PhD Candidate, School of Aerospace Engineering, lihy15@mails.tsinghua.edu.cn.} and Shiyu Chen\footnote{Graduate Student, School of Aerospace Engineering, chen-sy16@mails.tsinghua.edu.cn.}}
\affil{Tsinghua University, Beijing, 100084, China}
\author{Dario Izzo\footnote{Scientific Coordinator, European Space Research and Technology Center, Advanced Concepts Team.}}
\affil{European Space Research and Technology Center, 2201 AZ Noordwijk, The Netherlands}
\author{Hexi Baoyin\footnote{Professor, corresponding author, School of Aerospace Engineering; baoyin@tsinghua.edu.cn. Senior member AIAA.}}
\affil{Tsinghua University, Beijing, 100084, China}

\begin{document}

\maketitle

\begin{abstract}
In the design of multitarget interplanetary missions, there are always many options available, making it often impractical to optimize in detail each transfer trajectory in a preliminary search phase. Fast and accurate estimation methods for optimal transfers are thus of great value. In this paper, deep feed-forward neural networks are employed to estimate solutions to three types of optimization problems: the transfer time of time-optimal low-thrust transfers, fuel consumption of fuel-optimal low-thrust transfers, and the total $\Delta$\emph{v} of minimum-$\Delta$\emph{v} J$_2$-perturbed multi-impulse transfers. To generate the training data, low-thrust trajectories are optimized using the indirect method and J$_2$-perturbed multi-impulse trajectories are optimized using J$_2$ homotopy and particle swarm optimization. The hyper-parameters of our deep networks are searched by grid search, random search, and the tree-structured Parzen estimators approach. Results show that deep networks are capable of estimating the final mass or time of optimal transfers with extremely high accuracy; resulting into a mean relative error of less than 0.5\% for low-thrust transfers and less than 4\% for multi-impulse transfers. Our results are also compared with other off-the-shelf machine-learning algorithms and investigated with respect to their capability of predicting cases well outside of the training data. 
\end{abstract}

\section*{Nomenclature}

{\renewcommand\arraystretch{1.0}
\noindent\begin{longtable*}{@{}l @{\quad=\quad} l@{}}
$\emph{a}$ & semi-major axis \\
$\emph{e}$ & eccentricity \\
$\emph{i}$ & inclination \\
$\Omega$ & right ascension of ascending node\\
$\omega$ & perigee \\
$\mu$ & gravitational constant of central body \\
$u$ & engine thrust ratio \\
$\textbf{\emph{$\alpha$}}$ & unit vector of thrust direction \\
$g_0$ & standard value of gravitational acceleration \\
\emph{T$_{max}$} & maximum thrust magnitude \\
\emph{I$_{sp}$} & thruster specific impulse\\
\textbf{\emph{M}} & transformation matrix of modified equinoctial elements \\
\textbf{\emph{D}} & gravity vector \\
$\epsilon$ & homotopy parameter \\
$m$ & mass of spacecraft \\
$\lambda$ & co-state vector \\
$L$ & phase angle \\
\emph{n$_{layer}$}& number of layers\\
\emph{n$_{neuron}$}& number of neurons\\
\emph{f}& activation function\\
\emph{opt}& optimizer\\
$\eta$& learning rate\\
$c$& decay rate of learning rate\\
\emph{B}& batch size\\
\emph{F}& feature\\

\end{longtable*}}

\section{Introduction}
\lettrine{M}{ultitarget} interplanetary missions are a complex domain for optimization. In space missions to multiple targets, substantial benefits are gained since the average expense of exploring each target is lowered. Multitarget missions that are already being conducted are focused mostly on small celestial bodies in our solar system (except multi-gravity-assist missions). Small celestial bodies attract vast research interest because of their great significance in various aspects \cite{yang2017rapid,zeng2016solar,cheng2017asteroid,wie2017planetary}. NASA’s Near Earth Asteroid Rendezvous (NEAR) Shoemaker performed a flyby of the asteroid 253 Mathilde on the way to the asteroid 433 Eros \cite{prockter2002near}. In the previous plan for the NEAR mission, an ambitious plan called Small-Body Grand Tour, which aimed to achieve flybys of two comets and two asteroids over a 10-year period, was proposed \cite{farquhar1993extended}. Deep Space 1 \cite{rayman1999mission} and Dawn \cite{russell2012dawn} also achieved multitarget visits of small celestial bodies. The Global Trajectory Optimization Competitions (GTOCs), which are one of the most challenging events in space engineering, had several editions greatly focussed on multitarget small celestial body missions\footnote{\url{https://sophia.estec.esa.int/gtoc_portal/}}. In addition to multitarget missions to small celestial bodies, multitarget active debris removal (ADR) missions have also been avidly studied recently. The population of low-earth-orbit (LEO) debris has grown rapidly in the last two decades \cite{liou2008instability} and will continue to increase even if no new spacecraft is launched, as a consequence of the collision between existing debris, known as Kessler’s syndrome \cite{kessler2010kessler}. Many ADR methods are being studied, including contactless interaction \cite{li2018dynamics} and the debris engine\footnote{Lan, L., Li, J., and Baoyin, H., “Debris Engine: A Potential Thruster for Space Debris Removal,” 2015, \url{https://www.technologyreview.com/s/ 544156/junk-eating-rocket-engine-could-clear-space-debris/.}}. Optimization of multitarget ADR missions is a complex mathematical problem that can be seen as a complex variant to the more famous Travelling Salesman Problem \cite{izzo2015evolving} and was the subject of the 9th GTOC \cite{izzo2018kessler} and 8th China Trajectory Optimization Competition, of which the best solutions are presented in \cite{gtoc9JPL, li2017j2}.

The preliminary phases in the design of multitarget space missions often consider these problems as combinatorial-continuous \cite{yang2018optimization}. The combinatorial part is a sequence optimization problem in which many mission options are available, and the continuous part is a single-leg transfer trajectory optimization problem where two types of transfer models, low-thrust and multi-impulse, are mainly considered. The low-thrust trajectories can be optimized using indirect methods \cite{jiang2012practical,zhang2015low,jiang2017improving,chi2018power,chen2018multi}. The multi-impulse trajectories can be optimized using particle swarm optimization (PSO) \cite{pontani2010particle,li2017gtoc9} and genetic algorithms \cite{luo2007optimal}. Although these proposed methods are efficient, they are time-consuming when a large number of transfers are required to be optimized. Therefore, it is impractical to optimize each trajectory in the preliminary design of multitarget missions, and fast estimation methods are thus useful to approximate the optimal solutions of transfers between two targets as highlighted in \cite{hennes2016fast}.

The famous Edelbaum’s approximation gives the formulas of characteristic velocity requirements for multi-impulse and low-thrust transfers between quasi-circular orbits \cite{edelbaum1961propulsion}. In \cite{li2017j2}, the formulas are extended to J2-perturbed low-thrust problems. Edelbaum’s method is completely analytical, and therefore very suitable for a large-scale sequence search problem. However, for low-thrust optimization problems, these approximation formulas are derived under the assumption that the transfer is a many-revolution transfer. Many-revolution low-thrust transfers are impractical in deep-space multitarget missions because even one revolution will already take years (for example, 4 or 5 years in the main belt). Thus, transfers between two targets in multitarget deep space missions are usually less than one revolution. When the transfer time is short, Edelbaum’s methods become inaccurate, and estimation methods for short transfers are required. An approximate method is presented to estimate the velocity change of low-thrust transfers between co-planar, low-eccentricity orbits \cite{casalino2014approximate}, and is extended to transfers with small inclination change \cite{gatto2014fast}. These methods are based on Edelbaum’s approximation and are accurate for the evaluation of multitarget missions. In these methods, equations need to be iteratively solved and an additional correction is necessary for rendezvous transfers, making these methods somewhat complex to implement. In addition to the methods based on Edelbaum’s approximation, machine-learning methods \cite{hennes2016fast,mereta2017machine} are proposed to estimate the maximum final mass of the fuel-optimal rendezvous problem. In their research, after the features are well selected, different regressors are tested, and their performances are compared. Machine-learning methods are found to be vastly superior to other commonly used approximation methods, such as Lambert estimation. Estimation methods of multi-impulse optimization have rarely been studied in the past, but it is also time-consuming to optimize many-revolution J2-perturbed multi-impulse trajectories \cite{li2017gtoc9}. Estimation of multi-impulse optimization will gradually become necessary with the increment of studies in, for instance, multitarget ADR missions in low Earth orbit.

Deep neural networks (DNNs) have attracted the interest of many researchers and have had applications in various fields in recent years. DNNs consist of connected layers that are composed of neurons and are able to exhibit desired behavior after learning \cite{schmidhuber2015deep}. DNNs have their applications in the field of astrodynamics as well, including optimal control problems  \cite{sanchez2018real,izzo2018machine}, autonomous lunar landing \cite{furfaro2018deep}, and orbit prediction \cite{peng2018artificial}. In \cite{sanchez2018real}, DNNs are trained to learn the optimal control in four different cases of pinpoint landing. Their results, which show that the landings controlled by DNNs are close to optimal ones, allow for the design of on-board real-time optimal control. In \cite{izzo2018machine}, DNNs are trained to design the Earth-Mars orbital transfers of mass optimal control. In \cite{furfaro2018deep}, DNNs, i.e., convolutional neural networks and recurrent neural networks, are trained to perform an autonomous Moon landing using only raw images taken by on-board cameras, with the optimal control directly determined by the images only. In \cite{peng2018artificial}, DNNs are trained by historical orbit determination and prediction data to improve orbit prediction accuracy. It can be seen from the above research that DNNs can establish good connections between input and output, which is very difficult to obtain in an analytical way. The estimation of optimal transfers considered in this paper is also a problem for which a connection between the transfer parameters and the optimal solution must be established. Therefore, DNNs are employed in this paper to estimate optimal solutions of transfers.

Datasets that are necessary for the training of DNNs are generated by optimizing transfer trajectories. In this paper, three different cases, time optimal low-thrust, fuel optimal low-thrust, and minimum $\Delta$\emph{v} J$_{2}$-perturbed multi-impulse trajectories, are optimized. For the low-thrust optimization, the orbital states are expressed in modified equinoctial elements (\emph{MEEs}) because they are non-singular and are thus efficient and robust \cite{li2017j2}. Co-state normalization \cite{jiang2012practical} is utilized to improve the efficiency of the initial guess of co-states. The logarithmic homotopic method \cite{jiang2012practical} is applied to overcome the difficulties caused by bang-bang control in fuel optimal problems. For J$_{2}$-perturbed multi-impulse optimization, a many-revolution J$_{2}$-Lambert problem is solved using J$_{2}$ homotopy at first. In J$_{2}$ homotopy, the constant J$_{2}$ is multiplied by a homotopic parameter $\epsilon$ that varies from 0 to 1 iteratively to approach the J$_{2}$-perturbed problem from a two-body-problem perspective. PSO is then applied to optimize the total $\Delta$\emph{v}. Then, 100,000 trajectories are generated for each case, and are divided into a training dataset, validation dataset, and test dataset. The training dataset is used for DNNs to learn. The validation dataset is used for model selection. There are many hyper-parameters in DNN models, such as learning rate and batch size. Hyper-parameters are set prior to learning and not selected by the learning algorithm itself. A hyper-parameter search is necessary because the performance of a DNN is highly affected by hyper-parameters. In this paper, three techniques, grid search \cite{bengio2012practical}, random search \cite{bergstra2012random}, and the tree-structured Parzen estimators (TPE) approach \cite{bergstra2011algorithms}, are applied to search the optimal DNN model. The well-trained DNNs are capable of estimating optimal transfers with extremely high accuracy. DNNs also show excellent performance when compared with other machine-learning algorithms and investigated based on generalization capability. Finally, DNNs are applied to three mult-itarget missions and prove to be highly fast and accurate.

The rest of this paper is organized as follows. In Section II, trajectory optimization methods are presented. In Section III, the datasets are generated using a main-belt asteroid rendezvous mission and multitarget ADR mission as backgrounds. In Section IV, the architecture of DNNs and hyper-parameter search techniques are detailed. In Section V, DNNs are analyzed and applied to multitarget missions. Conclusions are drawn in Section VI.

\section{Trajectory Optimization}

\subsection{Low-Thrust Trajectory Optimization}

Modified equinoctial elements (\emph{MEEs}) are used here as the states \textbf{\emph{x}} because they are non-singular and thus are efficient and robust. Knowing the classical orbital elements (\emph{COEs}), the \emph{MEEs} can be expressed as:
\begin{equation}
\label{mee}
\begin{array}{l}
p = a\left( {1 - {e^2}} \right),\\
{e_x} = e\cos \left( {\Omega  + \omega } \right),\\
{e_y} = e\sin \left( {\Omega  + \omega } \right),\\
{h_x} = \tan \left( {{i \mathord{\left/
 {\vphantom {i 2}} \right.
 \kern-\nulldelimiterspace} 2}} \right),\cos \Omega \\
{h_y} = \tan \left( {{i \mathord{\left/
 {\vphantom {i 2}} \right.
 \kern-\nulldelimiterspace} 2}} \right),\sin \Omega \\
L = \Omega  + \omega  + f,
\end{array}
\end{equation}
where \emph{a} is the semi-major axis, \emph{e} the eccentricity, \emph{i} the inclination, $\Omega$ the right ascension of the ascending node, $\omega$ the perigee, and \emph{f} the true anomaly.

The dynamic equations of low-thrust propelled spacecraft in a two-body model are:
\begin{equation}
\label{dynamic_lowthrust}
\begin{array}{l}
\dot {\textbf{\emph{x}}} = {\textbf{\emph{M}}}\dfrac{{{T_{\max }}u}}{m}{\textbf{\emph{$\alpha$}}}  + {\textbf{\emph{D}}},\\
\dot m =  - \dfrac{{{T_{\max }}}}{{{I_{sp}}{g_0}}}u,
\end{array}
\end{equation}
where \textbf{\emph{x}} = [\emph{p}, \emph{e$_x$}, \emph{e$_y$}, \emph{h$_x$}, \emph{h$_y$}, \emph{L}], \emph{m} is the instantaneous mass of the spacecraft, $\mu$ the gravitational constant of the central body, u (0 $\le$ u $\le$ 1) the engine thrust ratio, \textbf{\emph{$\alpha$}} the unit vector of the thrust direction, g$_0$ = 9.80665 m/s$^2$ the standard value of gravitational acceleration, \emph{T$_{max}$} the maximum achievable thrust magnitude, and \emph{I$_{sp}$} the thruster specific impulse. \textbf{\emph{M}} is a transformation matrix and \textbf{\emph{D}} the gravity vector, the details of which can be found in \cite{gao2004low}.

For time-optimal problems (TOPs) and fuel-optimal problems (FOPs), the performance indexes take the forms:
\begin{equation}
\label{J}
\begin{array}{l}
J = {\lambda _0}\int_{{t_0}}^{{t_f}} {Ldt}, \\
{L_{{\rm{TOP}}}} = 1,\\
{L_{{\rm{FOP}}}} = \dfrac{{{T_{max}}}}{{{I_{sp}}{g_0}}}\left( {u - \varepsilon \ln \left[ {u\left( {1 - u} \right)} \right]} \right),
\end{array}
\end{equation}
%
%
%
%
where $\lambda_0$ is a positive factor that does not change the origin optimal problem and is necessary for the normalization of the co-state variables \cite{jiang2012practical}. For FOPs, the performance index uses the logarithmic homotopic method \cite{bertrand2002new} in which the fuel optimal index is achieved when $\epsilon$ becomes 0 from 1. Introducing the co-state vector \textbf{$\lambda$}(\emph{t}) = [\textbf{$\lambda$}$_{\textbf{\emph{x}}}$, \textbf{$\lambda$}$_m$, \textbf{$\lambda$}$_0$], the Hamiltonian is built as
\begin{equation}
\label{H}
H = \textbf{\emph{$\lambda$}} _{\textbf{\emph{x}}}^T \textbf{\emph{M}} \dfrac{{{T_{max }}u}}{m}\textbf{\emph{$\alpha$}}  + \textbf{\emph{$\lambda$}} _{\textbf{\emph{x}}}^T \textbf{\emph{D}} - {\lambda _m}\dfrac{{{T_{\max }}u}}{{{I_{sp}}{g_0}}} + {\lambda _0}L.
\end{equation}

According to Pontryagin’s minimum principle, the direction of the optimal thrust is determined as follows:
\begin{equation}
\label{alpha}
\textbf{\emph{$\alpha$}}  =  - \frac{{{\textbf{\emph{M}}^T}\textbf{\emph{$\lambda$}} _{\textbf{\emph{x}}}}}{{\left\| {{{\textbf{\emph{M}}^T}\textbf{\emph{$\lambda$}} _{\textbf{\emph{x}}}}} \right\|}}.
\end{equation}

The optimal thrust magnitude of TOPs is:
\begin{equation}
\label{utop}
\begin{array}{l}
{u_{{\rm{TOP}}}} = \left\{ {\begin{array}{*{20}{l}}
{0,{\rm{if}}{\rho _{{\rm{TOP}}}} > 0},\\
{1,{\rm{if}}{\rho _{{\rm{TOP}}}} < 0},\\
{[0,1],{\rm{if}}{\rho _{{\rm{TOP}}}} = 0},
\end{array}} \right.\\
{\rho _{{\rm{TOP}}}} =  - \dfrac{{{I_{{\rm{sp}}}}{g_0}\left\| {{\textbf{\emph{M}}^T}\textbf{\emph{$\lambda$}} _{\textbf{\emph{x}}}} \right\|}}{{{\lambda _0}m}} - \dfrac{{{\lambda _m}}}{{{\lambda _0}}},
\end{array}
\end{equation}
and the optimal thrust magnitude of FOPs is:
\begin{equation}
\label{ufop}
\begin{array}{l}
{u_{{\rm{FOP}}}} = \dfrac{{2\varepsilon }}{{{\rho _{{\rm{FOP}}}} + 2\varepsilon  + \sqrt {\rho _{{\rm{FOP}}}^2 + 4{\varepsilon ^2}} }},\\
{\rho _{{\rm{FOP}}}} = 1 - \dfrac{{{I_{{\rm{sp}}}}{g_0}\left\| {{\textbf{\emph{M}}^T}\textbf{\emph{$\lambda$}} _{\textbf{\emph{x}}}} \right\|}}{{{\lambda _0}m}} - \dfrac{{{\lambda _m}}}{{{\lambda _0}}}.
\end{array}
\end{equation}

Given the boundary condition, the optimal control problem can be transformed into a two-point boundary-value problem (TPBVP) and can be solved by shooting methods solving shooting equations. For rendezvous problems, in which the final time is optimized in TOPs and fixed in FOPs, the state variables of the spacecraft satisfy the following boundary conditions:
\begin{equation}
\label{bc}
\begin{array}{c}
\textbf{\emph{x}}({t_0}) = \textbf{\emph{x}}{_0},
m({t_0}) = {m_0},\\
\textbf{\emph{x}}({t_f}) = \textbf{\emph{x}}{_f}.
\end{array}
\end{equation}

The shooting functions of TOPs and FOPs can be deduced according to optimal control theory \cite{chi2018power,li2017j2} and are given as follows:

\begin{equation}
\label{shootT}
\textbf{\emph{$\Phi$}}{( {\textbf{\emph{$\lambda$}}({{t_0}})})_{{\rm{TOP}}}} = {[\textbf{\emph{x}}( {{t_f}}) - {\textbf{\emph{x}}_f}, {\lambda _m}( {{t_f}} ), ||{{\textbf{\emph{$\lambda$}}( {{t_0}} )}|| - 1}, H( {{t_f}}) - {\textbf{\emph{$\lambda$}} _x}( {{t_f}}) \cdot {\dot{\textbf{\emph{x}}}_f}]^T } = 0,
\end{equation}
\begin{equation}
\label{shootf}
\textbf{\emph{$\Phi$}}{( {\textbf{\emph{$\lambda$}}({{t_0}})})_{{\rm{FOP}}}} = {[\textbf{\emph{x}}( {{t_f}}) - {\textbf{\emph{x}}_f}, {\lambda _m}( {{t_f}} ), ||{{\textbf{\emph{$\lambda$}}( {{t_0}} )}|| - 1}]^T} = 0.
\end{equation}

Compared with FOPs, TOPs are relatively easy to solve due to the continuity and differentiability of the right-hand sides of the dynamical equations and Euler-Lagrange equations. A number of ideas have been advanced, also recently~\cite{jiang2012practical,jiang2017improving}, to improve the FOP convergence.

\subsection{J2-perturbed Multi-Impulse Trajectory Optimization}

For a higher computing efficiency of J2-perturbed propagation, the dynamic equations are set up using mean orbit elements (\emph{MOEs}):
\begin{equation}
\label{moe}
\begin{array}{l}
\dfrac{{{\rm{d}}{a_m}}}{{{\rm{d}}t}} = 0, \dfrac{{{\rm{d}}{e_m}}}{{{\rm{d}}t}} = 0, \dfrac{{{\rm{d}}{i_m}}}{{{\rm{d}}t}} = 0,\\
\dfrac{{{\rm{d}}{\Omega _m}}}{{{\rm{d}}t}} =  - \dfrac{3}{2}{J_2}{\left( {\dfrac{{{R_e}}}{p}} \right)^2}n\cos \left( {{i_m}} \right),\\
\dfrac{{{\rm{d}}{\omega _m}}}{{{\rm{d}}t}} = \dfrac{3}{4}{J_2}{\left( {\dfrac{{{R_e}}}{p}} \right)^2}n\left( {5{{\cos }^2}\left( {{i_m}} \right) - 1} \right),\\
\dfrac{{{\rm{d}}{M_m}}}{{{\rm{d}}t}} = n + \dfrac{3}{4}{J_2}{\left( {\dfrac{{{R_e}}}{p}} \right)^2}n\sqrt {1 - e_m^2} \left( {3{{\cos }^2}\left( {{i_m}} \right) - 1}, \right)
\end{array}
\end{equation}
where J$_2$=1.08262668e-3, $\emph{p}=a_m(1-e^2)$,$n=\sqrt[]{\mu/a_m^3}$, $\mu=3.9860044e14m^3/s^2$ is the gravitational constant of Earth, and the Earth radius R$_e$=6378137.0 m. The transformation between \emph{MOEs} and \emph{COEs} can be found in \cite{alfriend2009spacecraft}, and then the position \textbf{\emph{r}} and velocity \textbf{\emph{v}} can be acquired.
%
%
%
%
In the J$_2$-perturbed multi-impulse problems (JMPs) considered in this paper, the initial and final states, as well as the initial and final epochs \emph{t$_s$} and \emph{t$_N$}, are fixed. The entire trajectory of the spacecraft can be split into N+1 legs by N intermediate impulses. The position and velocity of the spacecraft at the epoch \emph{t$_i$} are denoted \textbf{\emph{r}}\emph{$_i$} and \textbf{\emph{v}}\emph{$_i$}$^\pm$, respectively, where the superscripts – and + represent the velocity before and after the impulse, respectively. Given the impulse vectors $\Delta$\textbf{\emph{v}}\emph{$_i$} and corresponding epochs \emph{t$_i$}, i=0,1,2,…,N-2, the first N legs can be propagated sequentially with Eq.~\eqref{moe}. The last leg is obtained by solving a J$_2$-Lambert problem using the shooting method to ensure that the rendezvous conditions are met.

J$_2$ homotopy is employed here to overcome the difficulties in the convergence of the shooting process caused by the long transfer duration. The constant J$_2$ is multiplied by a homotopic parameter $\epsilon\in$ [0,1] that varies from 0 to 1 iteratively to approach the J$_2$-perturbed problem from the perspective of a two-body problem. We denote the target position and velocity at the final epoch as \textbf{\emph{r}}\emph{$_f$} and \textbf{\emph{v}}\emph{$_f$}, respectively. For the last leg, the shooting uses \textbf{\emph{v}}\emph{$_{N-1}$}$^+$ as the shooting variable and $\Delta$\textbf{\emph{r}}=\textbf{\emph{r}}\emph{$_N$}- \textbf{\emph{r}}\emph{$_f$} as the shooting function. A two-body Lambert problem in which the spacecraft transfers from \textbf{\emph{r}}\emph{$_{N-1}$} to \textbf{\emph{r}}\emph{$_N$} for a flight time of \emph{t$_N$}-\emph{t$_{N-1}$} is solved first to obtain the initial guess of \textbf{\emph{v}}\emph{$_{N-1}$}$^+$, and then the homotopy process is implemented. The pseudo-algorithm is given in Algorithm \ref{alg1}. After the J$_2$-Lambert problem is solved, the last two impulses as well as the total $\Delta$\textbf{\emph{v}} can be calculated.

\begin{algorithm}
\caption{J$_2$ homotopy}
\label{alg1}
\begin{algorithmic}
\STATE {\textbf{\emph{v}}$_{initialguess}$ = solve(2bodyLambert)}.
\STATE {$\epsilon$ = 0}.
\WHILE {true}
\STATE {update($\epsilon$)};
\STATE {\textbf{\emph{v}}$_{N-1}^+$ = solve(HomotopicJ2Lambert($\epsilon$, \textbf{\emph{v}}$_{initialguess}$))};
\IF {$\epsilon$==1}
\STATE {break}
\ELSE
\STATE {\textbf{\emph{v}}$_{initialguess}$ = \textbf{\emph{v}}$_{N-1}^+$}
\ENDIF
\ENDWHILE
\end{algorithmic}
\end{algorithm}

The optimization of $\Delta$\textbf{\emph{v}}$_{total}$, which is determined by $\Delta$\textbf{\emph{v}}\emph{$_i$} and \emph{t$_i$}, is a typical parameter optimization problem, and swarm intelligence algorithms can find optimal solutions efficiently \cite{pontani2010particle,li2017gtoc9}. Thus, a PSO is applied to optimize $\Delta$\textbf{\emph{v}}$_{total}$ of the multi-impulse rendezvous problems in this paper. Considering a rendezvous process in which N+1 (N$\ge$2) impulses are employed, the unknown variables contain $\Delta$\textbf{\emph{v}}\emph{$_i$}, i=0,1,2,…,N-2, and \emph{t$_i$}, i=0,1,2,…,N-1, so the dimension of this case is n=(N-1)×3+N=4N-3. All the n unknown variables are denoted by the normalized position vector of a particle \textbf{\emph{x}}=[\emph{x$_1$}, \emph{x$_2$},…,\emph{x$_n$}], of which each component is in [0, 1]:
\begin{equation}
\label{pso1}
\begin{array}{l}
\Delta {\textbf{\emph{v}}_{\rm{i}}} = \Delta {v_{\rm{i}}}\left[ \begin{array}{c}
\cos \left( {{\alpha _{\rm{i}}}} \right)\cos \left( {{\beta _{\rm{i}}}} \right)\\
\cos \left( {{\alpha _{\rm{i}}}} \right)\sin \left( {{\beta _{\rm{i}}}} \right)\\
\sin \left( {{\alpha _{\rm{i}}}} \right)
\end{array} \right],\\
\Delta {v_{\rm{i}}} = \Delta {v_{\max }}{x_{3{\rm{i}} + 1}},\\
{\alpha _{\rm{i}}} = \pi \left( {{x_{3{\rm{i}} + 2}} - 0.5} \right),\\
{\beta _{\rm{i}}} = 2\pi {x_{3{\rm{i}} + 3}}\,\
i = 0,1, \ldots ,{\rm{N}} - 2,
\end{array}
\end{equation}

\begin{equation}
\label{pso2}
\begin{array}{l}
{t_0} = {t_{\rm{s}}} + {x_{3\left( {{\rm{N}} - 1} \right) + 1}}\left( {{t_{\rm{N}}} - {t_{\rm{s}}}} \right),\\
{t_{\rm{i}}} = {t_{{\rm{i}} - 1}} + {x_{3\left( {{\rm{N}} - 1} \right) + 1 + {\rm{i}}}}\left( {{t_{\rm{N}}} - {t_{{\rm{i}} - 1}}} \right),\\
i = 1,2, \ldots ,{\rm{N}} - 1.
\end{array}
\end{equation}

When N=1, which means only two impulses are employed, the dimension of \textbf{\emph{x}} is n=1 and the only unknown variable is t$_0$. The number of impulses (or N) is also unknown and it may differ for optimal solutions with different transfer conditions. To determine the optimal number of impulses, an enumeration for N is applied in this algorithm. Using PSO to obtain the optimal $\Delta$\textbf{\emph{v}}$_{total}$ for the cases where N=1, 2, 3, and 4 separately, the case with the smallest $\Delta$\textbf{\emph{v}}$_{total}$ is considered the optimal solution. Such an enumeration is implemented by parallel computing to improve the computational efficiency; that is, the four cases are computed concurrently through four threads.

\section{Dataset Generation}

\subsection{Dataset Generation of Low-Thrust Trajectories}

A main-belt asteroid rendezvous mission is considered as the test case for the low-thrust optimization problem. Asteroid orbital parameters are taken from the GTOC7 competition \cite{casalino2014problem} description. Altogether, 100,000 time-optimal trajectories and 100,000 fuel-optimal trajectories are generated. To make the generation process fast and highly efficient, asteroids with eccentricities larger than 0.1 are excluded. The process of low-thrust database generation is given in Algorithm \ref{alg2}.

\begin{algorithm}
\caption{Generation of Low-Thrust Trajectories}
\label{alg2}
\begin{algorithmic}
\WHILE {data number is less than 100,000}
\STATE {select two asteroids A1 and A2 randomly};
\IF {abs(\emph{a}$_{A1}$-\emph{a}$_{A2}$)>\emph{a}$_{MAX}$ or abs(\emph{i}$_{A1}$-\emph{i}$_{A2}$)>\emph{i}$_{MAX}$}
\STATE {continue};
\ENDIF
\STATE {give start time randomly}
\IF {abs(\emph{L}$_{A1}$-\emph{L}$_{A2}$)>\emph{L}$_{MAX}$}
\STATE {continue};
\ENDIF
\STATE {give initial mass randomly};
\STATE {solve TOP}
\IF {not converge}
\STATE {continue}
\ENDIF
\STATE {give the transfer time of FOP randomly}
\STATE {solve FOP}
\IF {not converge}
\STATE {continue}
\ENDIF
\STATE {add these data to database}
\ENDWHILE
\end{algorithmic}
\end{algorithm}

A pre-screening procedure is adopted to select asteroid pairs. The \emph{a}$_{MAX}$ is 0.3 AU, \emph{i}$_{MAX}$ is 3 deg, and \emph{L}$_{MAX}$ is 30 deg. This procedure ensures that the two selected asteroids are close enough to each other at the given epoch to improve the convergence rate of low-thrust optimal control problems. The start time is randomly given from [56800, 59800] MJD, and the initial mass is randomly given from [1000, 2000] kg. If the TOP converges, the FOP is then solved with the transfer time transfer randomly given from [1.2, 1.5]×\emph{t}$_{top}$. If also the FOP converges, then the following information will be saved to the database:

\emph{A1, A2, start time, initial mass, transfer time of FOP, optimal time of TOP, optimal fuel of FOP.}

There are many local optima in complex optimization problems. In the indirect method of solving low-thrust optimization problems, different initial guesses of co-states can lead to different local optima. Therefore, each TOP and FOP is solved 10 times and the best result is saved as the optimal value.

\subsection{Dataset Generation of Multi-Impulse Trajectories}

A multitarget ADR mission is considered as the test case for the multi-impulse optimization problem. Debris orbital parameters are from the GTOC9 competition \cite{izzo2018kessler}. Altogether, 100,000 J2-perturbed multi-impulse trajectories are generated. The process of multi-impulse database generation is given in Algorithm \ref{alg3}.

\begin{algorithm}
\caption{Generation of J$_2$-Perturbed Multi-Impulse Trajectories}
\label{alg3}
\begin{algorithmic}
\WHILE {data number is less than 100,000}
\STATE {select two debris D1 and D2 and give start time randomly};
\IF {estimate$\Delta$\emph{v}>$\Delta$\emph{v}$_{MAX}$}
\STATE {continue};
\ENDIF
\STATE {given transfer time randomly}
\STATE {solve multi-impulse optimization using PSO};
\STATE {add these data to database}
\ENDWHILE
\end{algorithmic}
\end{algorithm}

A $\Delta$\emph{v} estimation calculated by Eq. (24) \cite{li2017gtoc9} is used as pre-screening. If this $\Delta$\emph{v} is larger than $\Delta$\emph{v}$_{max}$= 500 m/s, than this debris pair will be abandoned:

\begin{equation}
\label{dv}
\Delta v = \left( {\sqrt {{{\left( {\dfrac{{\Delta a}}{{2{a_1}}}} \right)}^2} + {{\left( {\dfrac{{\Delta e}}{2}} \right)}^2} + {{\left( {\Delta i} \right)}^2}}  + \Delta \Omega \sin {i_1}} \right)\sqrt {\dfrac{{2\mu }}{{{a_1}}}}.
\end{equation}

The start time is randomly given from [56,800, 59,800] MJD, and the transfer time is randomly given from [5, 30] days. The following information will be saved to the database:

\emph{A1, A2, start time, transfer time, optimal $\Delta$\emph{v}.}

The PSO will run three times for each problem and the best result is saved as the optimal value.

\section{Neural Architecture}

\subsection{Architecture of DNNs}

\begin{figure}[hbt!]
\centering
\includegraphics[width=.7\textwidth]{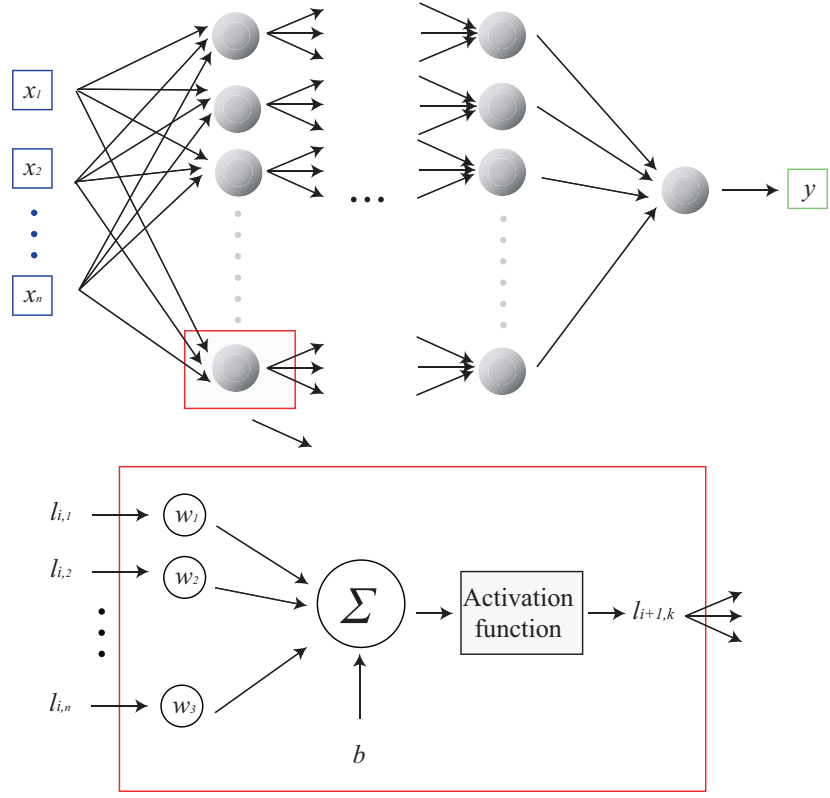}
\caption{Feed-forward neural network.}
\label{figDNN}
\end{figure}

The structure of the network considered in this paper is illustrated in Fig.\ref{figDNN} and is essentially a feed-forward neural network with multiple hidden layers. The structure is determined by the number of layers \emph{n}$_{layer}$ and the number of neurons \emph{n}$_{neuron}$ at each hidden layer. The data will flow through these layers and finally produce the output. At each layer, denoting the input \emph{l$_i$}, the output \emph{l$_{i+1}$} is calculated as follows:

\begin{equation}
\label{layer}
{\textbf{\emph{{l}$_{i + 1}$}}} = \textbf{\emph{f}}( {\textbf{\emph{w}}{\textbf{\emph{l}}_i} + \textbf{\emph{b}}}),
\end{equation}
where \textbf{\emph{w}} is the weight matrix, \textbf{\emph{b}} the bias vector, and \textbf{\emph{f}} a nonlinear function named the activation function. This mechanism is also shown in Fig.\ref{figDNN}. There are three most commonly used activation functions for the hidden layers: the sigmoid function expressed in Eq.\eqref{sig}, the hyperbolic tangent (tanh) function expressed in Eq.\eqref{tanh}, and the rectified linear (ReLu) function expressed in Eq.\eqref{relu}. For the output layer, a linear function is selected:

\begin{equation}
\label{sig}
f_{sig}(x) = \dfrac{1}{{1 + {e^{ - x}}}},
\end{equation}

\begin{equation}
\label{tanh}
f_{tanh}(x) = \dfrac{{{e^x} - {e^{ - x}}}}{{{e^x} + {e^{ - x}}}},
\end{equation}

\begin{equation}
\label{relu}
f_{relu}(x) = \max ({0,x}).
\end{equation}

%
%
%
%

The training process consists of adjusting the value of the parameters of each layer to minimize the loss function:

\begin{equation}
\label{loss}
E = \dfrac{1}{n}\sum\limits_{i = 1}^n {{{\left( {\hat y - y} \right)}^2}},
\end{equation}
where y is the actual value,  $\hat y$ the estimated value, and n the number of data in this training iteration.  To evaluate the performance of DNNs, the error function is given as:
\begin{equation}
\label{error}
\begin{array}{l}
{\rm{error}} = \dfrac{1}{n}\sum\limits_{i = 1}^n {\dfrac{{\left| {\hat y - y} \right|}}{y}}, \\
{\rm{accuracy}} = 1 - {\rm{error}}.
\end{array}
\end{equation}

Gradient descent (GD) algorithms are the state-of-art in training the parameters, e.g.,

\begin{equation}
\label{GD}
\textbf{\emph{w}}' = \textbf{\emph{w}} - \eta \dfrac{{\partial E}}{{\partial \textbf{\emph{w}}}},
\end{equation}
where $\eta$ is the learning rate. Some modified gradient descent algorithms are also very effective, such as momentum gradient descent (MGD) \cite{sutskever2013importance} and Adam gradient descent (AGD) \cite{kingma2014adam}.

A most important hyper-parameter is the learning rate $\eta$. With a small learning rate, the DNN may take a very long time to converge, while with a large learning rate the DNN may oscillate around the optimal solution. Therefore, it is a good strategy to start from a large learning rate and decrease the learning rate gradually. Two decay schedules are applied in this paper: one is exponential decay (ED),

\begin{equation}
\label{ED}
\eta ' = \eta {c ^t},
\end{equation}
and the other is natural exponential decay (NED),
\begin{equation}
\label{NED}
\eta ' = \eta {e^{ - c t}},
\end{equation}

where \emph{c} is the decay rate and
\begin{equation}
\label{tstep}
t = \frac{{{global\_step}}}{{{decay\_step}}}.
\end{equation}

In training, the batch size B is also a hyper-parameter that must be chosen. A batch of data is put into the DNN to train it and then the weights are updated. This hyper-parameter effects training time more \cite{bengio2012practical}. Early stopping is utilized here, and the training will stop when there is no improvement in the last N epochs.

The form of input, which is also known as features, is also an important factor that affects the performance of DNNs. For the problem considered in this paper, the input is a set of transfer parameters. The orbit states can be expressed in many different forms \cite{hintz2008survey}. Three forms are considered in this paper: Cartesian form \emph{rv} (\emph{x}, \emph{y}, \emph{z}, \emph{v$_x$}, \emph{v$_y$}, \emph{v$_z$}), classical orbit elements \emph{COEs} (\emph{a}, \emph{e}, \emph{i}, \emph{$\Omega$}, \emph{$\omega$}, \emph{f}), and modified equinoctial elements \emph{MEEs} (\emph{p}, \emph{e$_x$}, \emph{e$_y$}, \emph{h$_x$}, \emph{h$_y$}, \emph{L}), and the corresponding input forms are donated F$_{rv}$, F$_{coe}$, and F$_{mee}$, respectively.

In addition to the basic transfer parameters, other information can also be added into the input. The differences of two integration constants, angular momentum h and mechanical energy E in Eq. \eqref{Eh}, between initial and final orbits, can be used as additional features:
\begin{equation}
\label{Eh}
\begin{array}{l}
\textbf{\emph{h}} = \textbf{\emph{r}} \times \textbf{\emph{v}},\\
E =  - \dfrac{\mu }{{2a}}.
\end{array}
\end{equation}

For the fuel-optimal, low-thrust problem, a two-body Lambert solution is used as an additional feature. For the J2-perturbed multi-impulse problem, Eq.\eqref{dv} is used as an additional feature. The input forms with additional features are donated F\emph{a}$_{rv}$, F\emph{a}$_{coe}$, and F\emph{a}$_{mee}$. All of the inputs are scaled to be between -1 and 1.
The actual transfer time (day) of TOPs is multiplied by a factor of 0.0172 as the actual output, the actual final mass (kg) of FOPs is multiplied by a factor of 0.004 as the actual output, and the actual $\Delta\emph{v}_{total}$ (m/s) of JMPs is multiplied by a factor of 0.01 as the actual output.

\subsection{Selection of DNN models}
The performance of a DNN is highly affected by hyper-parameters, which are set prior to the learning and not selected by the learning algorithm itself. Many model choices are implemented by manual search, which gives the impression of neural network training as an art \cite{bengio2012practical}. There are two drawbacks of manual search. One is that it is not easily reproducible and automated, and the other is that a human may never find the optimal hyper-parameters. Therefore, three techniques, grid search \cite{bengio2012practical}, random search \cite{bergstra2012random}, and the tree-structured Parzen estimators (TPE) approach \cite{bergstra2011algorithms}, are applied in this paper to select the most proper DNN model.

\emph{Grid Search}: Grid search is an automated search process executed after defining a search space. However, an exponential explosion may occur when the number of hyper-parameters increases. One must therefore define the search space carefully due to the limitation of computing capability.

\emph{Random Search}: Random search is highly recommended by \cite{bergstra2012random}, as they found random sampling very efficient. Random search allows one to explore more hyper-parameters and more values of one hyper-parameter, compared to grid search. The enlarged search space increase the possibility of finding the optimal hyper-parameters, while random sampling does not increase the computational budget. Random search is also simple to implement.

\emph{Tree-structured Parzen Estimators Approach}: TPE is a sequential model-based optimization (SMBO) algorithm. In grid search and random search, the next step is completely independent of previous steps. SMBO algorithms will consider previous results and decide what to try next.  TPE uses random forest regression to provide the prediction of which region of hyper-parameters is optimal. The details of TPE implementation can be found in \cite{bergstra2011algorithms}.

The hyper-parameters that must be searched are number of layers \emph{n$_{layer}$}, number of neurons at each hidden layer \emph{n$_{neuron}$}, activation function \emph{f}, batch size \emph{B}, optimizer \emph{opt}, initial learning rate $\eta$, decay model \emph{dm}, decay step \emph{ds}, and decay rate \emph{c}. In addition to the hyper-parameters, six feature forms, F$_{rv}$, F$_{coe}$,  F$_{mee}$, F\emph{a}$_{rv}$, F\emph{a}$_{coe}$, and F\emph{a}$_{mee}$, also must be searched. The search space of grid search is listed in Table \ref{tabGS} and that of random search and TPE in Table \ref{tabRSTPE}. There are altogether 96 trails in grid search; therefore, 100 trails are set in random search and TPE, so the three search methods are comparable. Even though their computing budgets are the same, random search and TPE can explore a larger parameter space. These parameters are uniformly random, except that the learning rate is uniform in the log-domain [40]. The other default hyper-parameters of grid search are \emph{f} = relu, \emph{opt} = GD, \emph{dm} = ED, \emph{ds} = 1000, and \emph{c} = 0.9. The step for early stopping is set to 50.

\begin{table}[hbt!]
\caption{\label{tabGS} Search space of grid search}
\centering
\begin{tabular}{lcccccc}
\hline
Hyper-parameters& Search space\\\hline
\emph{n$_{layer}$}& 2, 3, 4, 5\\
\emph{n$_{neuron}$}& 32, 128\\
$\eta$& 0.01, 0.001\\
\emph{B}& 200, 800\\
\emph{F}& F$_{rv}$, F$_{coe}$,  F$_{mee}$\\
\hline
\end{tabular}
\end{table}

\begin{table}[hbt!]
\caption{\label{tabRSTPE} Search space of random search and TPE}
\centering
\begin{tabular}{lcccccc}
\hline
Hyper-parameters& Search space\\\hline
\emph{n$_{layer}$}& 2, 3, 4, 5\\
\emph{n$_{neuron}$}& 32, 128\\
\emph{f}& sigmoid, tanh, relu\\
\emph{opt}& GD, MGD, AGD\\
$\eta$& 0.1 - 0.0001\\
\emph{dm}& ED, NED\\
\emph{ds}& 500 - 1500\\
$c$& 0.8 - 1\\
\emph{B}& 100 - 1000\\
\emph{F}& F$_{rv}$, F$_{coe}$,  F$_{mee}$, F\emph{a}$_{rv}$, F\emph{a}$_{coe}$, F\emph{a}$_{mee}$\\
\hline
\end{tabular}
\end{table}

The dataset is divided into 80,000 training datasets, 10,000 validation datasets, and 10,000 test datasets. The training dataset is used to train the network, which is inside the learning algorithm. The validation dataset is used for DNN model selection, which is outside the learning algorithm. The test dataset is finally used to evaluate the performance of selected model. Many take validation and test datasets as one dataset, but when a model-selection process is considered the validation and test datasets must be separated.

\section{Simulations}

\subsection{Dataset Analysis}

The basic dataset information is listed in Tables \ref{tab_dataset_lowthrust} and \ref{tab_dataset_impulse}. 
In Table \ref{tab_dataset_lowthrust}, the phase angle \emph{L} is calculated at the start time, and \emph{t}$_{transfer}$ is the transfer time of FOPs.
The maximum value of orbital differences in the low-thrust dataset is the same as the pre-screening settings. In Table \ref{tab_dataset_impulse}, the $\Delta\Omega$ is the orbital difference between start debris and target debris at the rendezvous time. Note that the minimum transfer time in Table \ref{tab_dataset_impulse} is less than 1 d, which is the lower boundary when generating the dataset, because the transfer time is the duration between rendezvous and first impulse time, which is an optimization variable.

\begin{table}[hbt!]
\caption{\label{tab_dataset_lowthrust} Basic information of low-thrust dataset}
\centering
\begin{tabular}{lcccccccc}
\hline
&\emph{m}$_s$ (kg)&|$\Delta$\emph{a}| (AU)&|$\Delta$\emph{e}|&$|\Delta$\emph{i}| (deg)&|$\Delta$\emph{L}| (deg)&Optimal time  (d)&\emph{t}$_{transfer}$&Optimal fuel (kg)\\\hline
mean&1490.2&0.1041&0.0262&1.31&14.85&533.76&719.60&270.15\\
max&2000&0.3&0.0971&3&30&1350.16&1791.60&873.76\\
min&1000&1e-6&1e-6&7e-5&8.4e-5&78.63&97.23&37.01\\
\hline
\end{tabular}
\end{table}

\begin{table}[hbt!]
\caption{\label{tab_dataset_impulse} Basic information of multi-impulse dataset}
\centering
\begin{tabular}{lcccccc}
\hline
&|$\Delta$\emph{a}| (km)&|$\Delta$\emph{e}|&$|\Delta$\emph{i}| (deg)&|$\Delta\Omega$| (deg)&\emph{t}$_{transfer}$&Optimal $\Delta$\emph{v} (m/s)\\\hline
mean&66.43&0.0055&0.76&1.48&7.18&119.51\\
max&277.88&0.0192&3.76&3.80&24.94&508.32\\
min&0.0079&3e-6&3.4e-4&5.3e-5&0.48&3.11\\
\hline
\end{tabular}
\end{table}

\begin{figure}[H]
\includegraphics[width=1.2\textwidth]{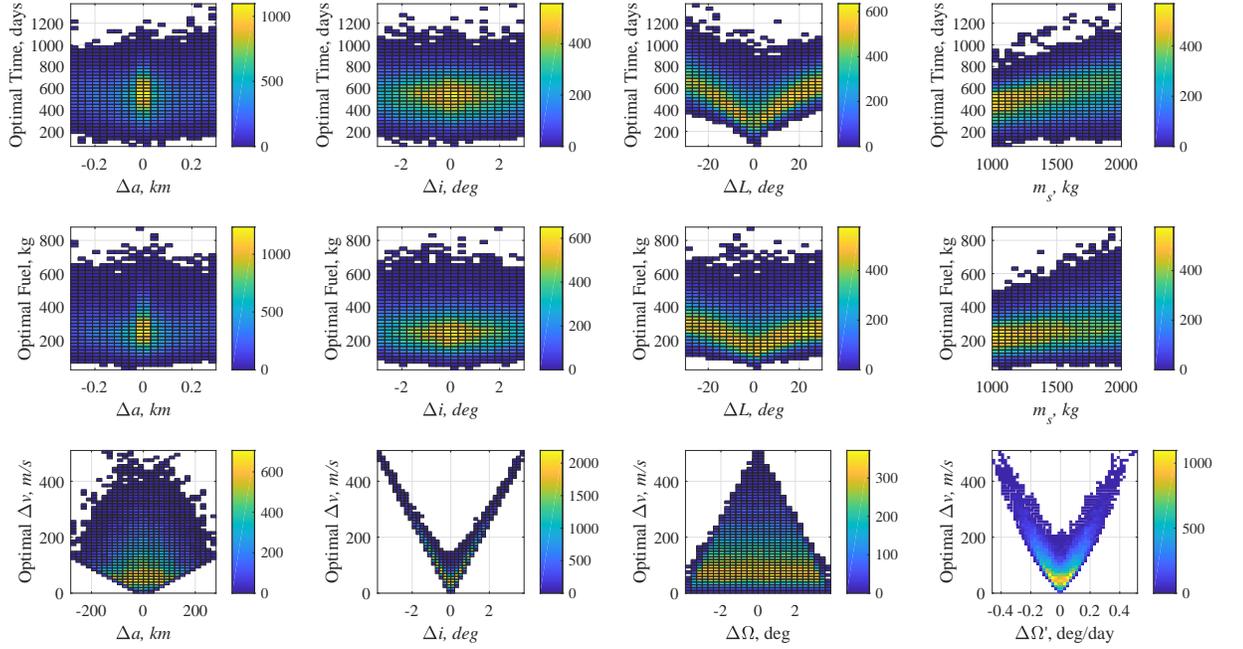}
\caption{Visualization of dataset content.}
\label{figDataset}
\end{figure}

The distribution of datasets is shown in Fig. \ref{figDataset}, where the upper, middle, and lower rows are TOPs, FOPs, and JMPs, respectively. The color bar indicates the amount.
In the last sub-figure of JMPs, $\Delta\Omega'$ is the difference of $\Omega$ precession rates between start debris and target debris.
For low-thrust transfers, it can be found that the phase angle has an apparent relevance to the optimal value. Transfers are concentrated in the center (small orbital differences) of the $\Delta\emph{a}$ and $\Delta\emph{i}$ sub-figures, and in the left (light in mass) of \emph{m$_s$} sub-figure because it is easy to converge for these transfers.
For J$_2$-perturbed multi-impulse transfers, a strong relevance can be found between $\Delta\emph{i}$ and the optimal value because the outer-plane maneuver costs $\Delta\emph{v}$ significantly. $\Delta\Omega'$ shows a similar relevance because $\Delta\Omega'$ is strongly determined by \emph{i}, which can be seen in Eq. \eqref{moe}. The $\Delta\Omega$ shows a nearly opposite tendency, although it is the outer-plane change. This is because the precession rate $\Omega'$caused by J$_2$ perturbation can be used to adjust the orbital plane instead of using large $\Delta\emph{v}$. The transfer with large $\Delta\Omega$ and small $\Delta\emph{i}$ can pass the pre-screening described in Algorithm \ref{alg3}, and the actual $\Delta\emph{v}$ is small because of the small $\Delta\emph{i}$.

\subsection{DNN Results}

\begin{figure}[H]
\centering
\includegraphics[width=0.9\textwidth]{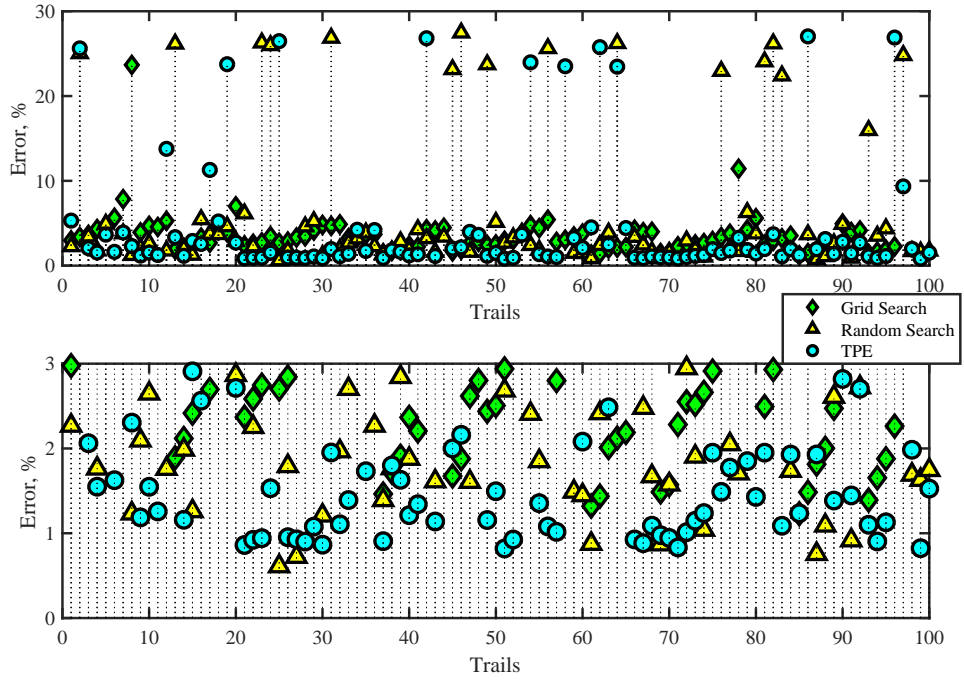}
\caption{Search results in time-optimal low-thrust problems.}
\label{figDatatop}
\end{figure}

\begin{figure}[H]
\centering
\includegraphics[width=0.9\textwidth]{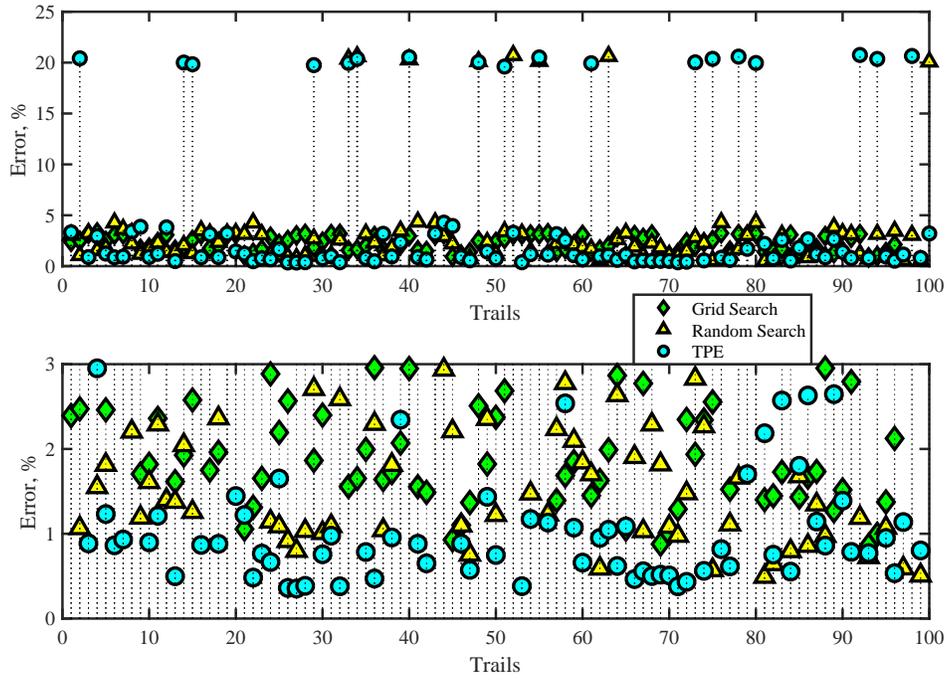}
\caption{Search results in fuel-optimal low-thrust problems.}
\label{figDatafop}
\end{figure}

\begin{figure}[H]
\centering
\includegraphics[width=0.9\textwidth]{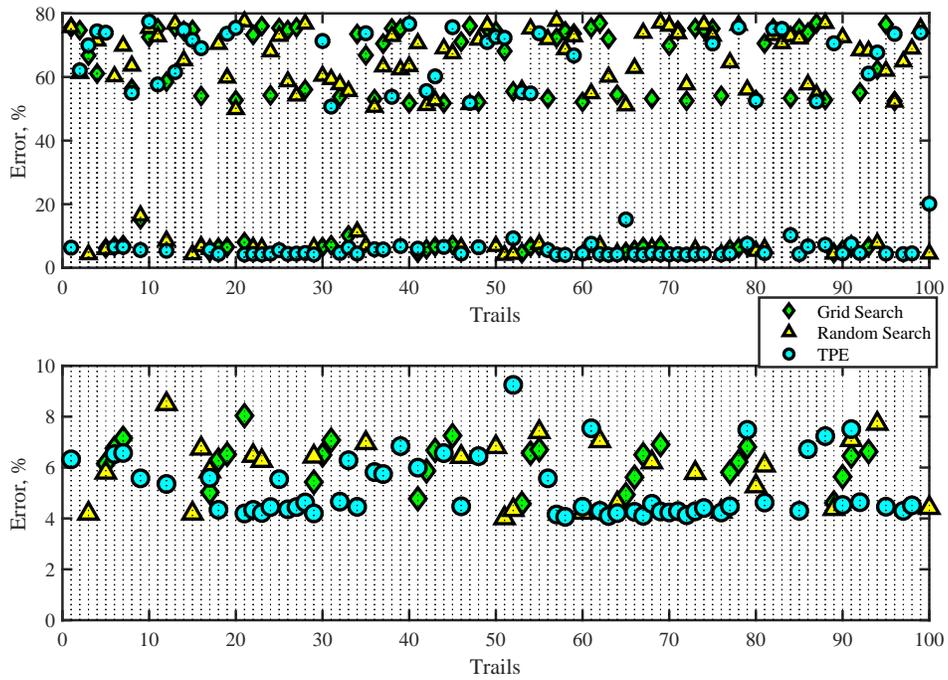}
\caption{Search results in J$_2$-perturbed multi-impulse problems.}
\label{figDataimp}
\end{figure}

Three hyper-parameter search techniques, grid search, random search, and TPE, are implemented to select the optimal DNN model. 
The results of TOPs, FOPs, and JMPs are shown in Figs. \ref{figDatatop}, \ref{figDatafop}, and \ref{figDataimp}, respectively.
In these figures, the x axis is the trails and the y axis the error on the validation dataset. It should be noted that the DNN model is trained on the training dataset and selected according to its performance on the validation dataset. Dots with three different colors and shapes represent the three search techniques. The lower sub-figure is the zoom of the bottom of the upper sub-figure.
For TOPs, the best DNN model is found by random search, and the corresponding error is 0.61\%. For FOPs, the best DNN model is found by TPE, and the corresponding error is 0.35\%. For JMPs, the best DNN model is found by random search, and the corresponding error is 4.01\%, while the error of the best model found by TPE is 4.06\%. Random search and TPE have similar performance, but both are superior to grid search. Random search is recommended because it is highly effective and easy to implement. 

Overall, the errors of JMPs are larger than those of TOPs and FOPs. This is because the optimality of low-trust datasets is better than that of the multi-impulse dataset. Among the methods used to optimize low-thrust trajectories, the indirect method utilized in this paper has the advantage of optimality. Although there are many local optima, the approach of optimizing the trajectory 10 times and choosing the best one can make the solution close to the global optimal solution. PSO as an stochastic optimization algorithm does not have an optimality as good as that of the indirect method. Therefore, the sub-optimal solutions in the multitarget dataset make the network difficult to learn.

Some errors are larger than 10\% because their learning parameters are not set properly. However, in Fig. \ref{figDataimp}, errors are larger than 60\% in some trails. By investigating their hyper-parameters, it is found that in most of these trails the features are F$_{rv}$ or F\emph{a}$_{rv}$. This is because, in long-duration multi-impulse transfers, the phase has very little effect on the optimal $\Delta\emph{v}$, but the phase does determine the position and velocity. Even if the spacecraft starts the transfer from two different position along the same orbit, which means the position and velocity are quite different, the optimal $\Delta\emph{v}$ can be nearly the same. Therefore, F$_{rv}$ and F\emph{a}$_{rv}$ are not suitable features for long-duration transfers.

The best DNN models selected by the search algorithms are listed in Table \ref{tab_hyper}. Some general guidance in choosing DNN models can be gleaned from this table, i.e., they are all quite deep and large networks that have several of the same or similar parameters.
The activation functions of the three problems are sigmoid, the optimizers are AGD, the initial learning rates are all 1e-3, and the decay models are NED. It should be noted that some of these hyper-parameters may never be selected by manual and grid search. 

\begin{table}[H]
\caption{\label{tab_hyper} Best DNN models.}
\centering
\begin{tabular}{lcccccccccc}
\hline
&\emph{n$_{layer}$}&\emph{n$_{neuron}$}&\emph{f}&\emph{opt}&$\eta$&\emph{dm}&\emph{ds}&$\lambda$&\emph{B}&\emph{F}\\\hline
TOP&4&218&sigmoid&AGD&1.99e-3&NED&700&0.85&250&F\emph{a}$_{rv}$\\
FOP&5&170&sigmoid&AGD&7.06e-3&NED&700&0.85&1000&F\emph{a}$_{ee}$\\
JMP&4&210&sigmoid&AGD&3.41e-3&NED&1000&0.83&850&F$_{coe}$\\
\hline
\end{tabular}
\end{table}

%
%
%
%

\begin{figure}[H]
\centering
\includegraphics[width=0.9\textwidth]{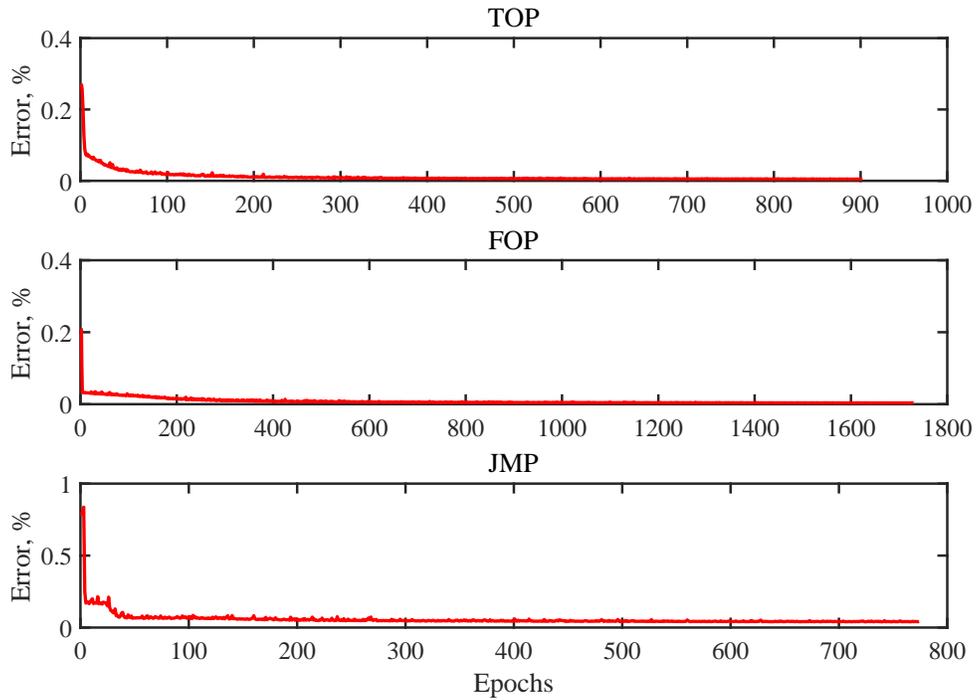}
\caption{Errors of test datasets in training process using selected DNN models.}
\label{figlearn}
\end{figure}

The performance of the selected DNN models must be checked on the test dataset. Since the initialization of weights and bias are random, each complete training process can lead to different results. Therefore, the DNNs are trained three times and that with the best performance on the test dataset will be saved for the following analysis. The training histories of the best DNNs of the three problems are shown in Fig. \ref{figlearn}. Up until now, well-trained DNNs have been created. In the following, accuracy will be used to present the performance of the DNNs. The accuracy on test datasets are listed in Table \ref{tab_learn}. For low-thrust problems, the accuracy of the DNNs is higher than 99.5\%, and for multi-impulse problems it is higher than 96\%. These results demonstrate that a well-trained DNN is capable of estimating optimal transfers with extremely high accuracy.

\begin{table}[H]
\caption{\label{tab_learn} Accuracy of test datasets.}
\centering
\begin{tabular}{lccc}
\hline
&1&2&3\\\hline
TOP&\textbf{99.54}\%&99.43\%&99.15\%\\
FOP&99.64\%&99.65\%&\textbf{99.66}\%\\
JMP&95.77\%&\textbf{96.06}\%&95.94\%\\
\hline
\end{tabular}
\end{table}

DNNs are compared with other machine-learning algorithms, e.g., those used in \cite{hennes2016fast,mereta2017machine}. 
A random search is utilized to tune these algorithms and the search space of the hyper-parameters is shown in Table \ref{tab_hpml}  \cite{scikit-learn}.
Features are the same as those in the best DNN models. The algorithms are trained on training datasets and their accuracy is tested on test datasets. For each case, 100 trails are implemented and the results are listed in Table \ref{tab_ml}. It is shown that DNNs are superior to all of the other machine-learning algorithms. However, these machine-learning algorithms also have excellent performance. In \cite{hennes2016fast,mereta2017machine}, estimated optimal transfers were already proved to have been improved significantly over those of the Lambert model, which is most commonly used. Our work confirms the effectiveness of machine learning, especially deep learning, in estimating optimal transfers.
%
%
%
%
\newcommand{\tabincell}[2]{\begin{tabular}{@{}#1@{}}#2\end{tabular}} 

\begin{table}[H]
\caption{\label{tab_hpml} Search space of other machine-learning algorithms.}
\centering
\begin{tabular}{lccc}
\hline
Algorithms&Search space\\\hline
RandomForest&\tabincell{c}{n\_estimators:	[10, 100], max\_depth: [2,100]\&None, \\min\_samples\_split: [2, 10], min\_samples\_leaf: [1, 4], bootstrap: [True, False]} \\
Bagging&n\_estimators: [10, 100], max\_samples: [0.2, 0.8], bootstrap: [True, False]\\
AdaBoost&n\_estimators: [10, 100], learning\_rate: [0.01, 2], loss: [linear, square]\\
ExtraTrees&\tabincell{c}{n\_estimators: [10, 100], criterion: [mse, mas], max\_depth: [2,100]\&None, \\min\_samples\_split: [2, 10], min\_samples\_leaf: [1, 4], bootstrap: [True, False]} \\
GradientBoosting&\tabincell{c}{n\_estimators: [10, 100],  loss: [ls, lad, huber], max\_depth: [2,100],\\ min\_samples\_split: [2, 10], min\_samples\_leaf: [1, 4],} \\
DecisionTree&\tabincell{c}{ splitter: [best, random], max\_depth: [2,100]\&None,\\ min\_samples\_split: [2, 10], min\_samples\_leaf: [1, 4], } \\
ExtraTree&\tabincell{c}{ splitter: [best, random], max\_depth: [2,100]\&None,\\ min\_samples\_split: [2, 10], min\_samples\_leaf: [1, 4], } \\
\hline
\end{tabular}
\end{table}

\begin{table}[H]
\caption{\label{tab_ml} Comparison with other machine-learning algorithms.}
\centering
\begin{tabular}{lccc}
\hline
&TOP&FOP&JMP\\\hline
RandomForest&94.75\%&98.11\%&95.10\%\\
Bagging&94.94\%&98.12\%&95.18\%\\
AdaBoost&86.37\%&96.33\%&82.92\%\\
ExtraTrees&95.50\%&98.29\%&95.22\%\\
GradientBoosting&95.15\%&98.16\%&95.29\%\\
DecisionTree&92.15\%&97.56\%&93.75\%\\
ExtraTree&91.98\%&97.41\%&93.68\%\\
DNN&\textbf{99.54}\%&\textbf{99.66}\%&\textbf{96.02}\%\\
\hline
\end{tabular}
\end{table}

To further investigate the performance of well-trained DNNs, generalization capability is studied by using the transfers in which the transfer parameters are outside the boundary of which values are given in Tables \ref{tab_dataset_lowthrust} and \ref{tab_dataset_impulse}.
For low-thrust problems, these transfers are generated by a similar procedure described in Algorithm \ref{alg2}, except that the pre-screening is different. In the pre-screening used to study generalization capability, the transfer will be saved if it is 1-1.2 times the boundary value. For instance, to study the generalization capability in terms of |$\Delta$\emph{a}|, the transfer will be saved only when |$\Delta$\emph{a}|$_{max}$ < |$\Delta$\emph{a}| < 1.2|$\Delta$\emph{a}|$_{max}$. To study the generalization capability in terms of \emph{m}$_s$, the initial mass is given from [2,000, 2,400] kg randomly.
When studying |$\Delta$\emph{e}|, asteroids with eccentricities larger than 0.1 are not excluded as before; otherwise, it will be very difficult for |$\Delta$\emph{e}| to be in the range 1-1.2 |$\Delta$\emph{e}|$_{max}$.
The initial mass \emph{m}$_s$, |$\Delta$\emph{a}|, |$\Delta$\emph{e}|, |$\Delta$\emph{i}|, and |$\Delta$\emph{L}| are studied separately, and for each situation 100 transfers are calculated.
The generalization capability results of low-thrust problems are given in Table \ref{tab_gc_topfop}. The accuracy is still higher than 99\% in terms of |$\Delta$\emph{a}|, |$\Delta$\emph{i}|, and |$\Delta$\emph{L}|. In |$\Delta$\emph{e}|, the accuracy is 97\%, because including asteroids with eccentricities larger than 0.1 also brings in the change of the original problems. 

For multi-impulse problems, |$\Delta$\emph{a}|, |$\Delta$\emph{e}|, |$\Delta$\emph{i}|, and |$\Delta\Omega$| are studied separately. For the study of |$\Delta$\emph{i}| and |$\Delta\Omega$|, transfers are generated by a similar procedure as described in Algorithm \ref{alg3}. In addition to |$\Delta$|$_{max}$ < |$\Delta$| < 1.2|$\Delta$|$_{max}$ being necessary, $\Delta\emph{v}_{max}$ is changed to 700 m/s. Since |$\Delta$\emph{a}|$_{max}$ and |$\Delta$\emph{e}|$_{max}$ almost reach the limitation of all the 123 debris, the approach to increasing |$\Delta$\emph{a}| or |$\Delta$\emph{e}| is to select one transfer from a dataset randomly and change the \emph{COE} of the target orbit to make |$\Delta$| in the range 1-1.2 times |$\Delta$|$_{max}$. The generalization capability results of multi-impulse problems are given in Table \ref{tab_gc_imp}, and demonstrate that well-trained DNNs have an excellent generalization capability.

\begin{table}[H]
\caption{\label{tab_gc_topfop} Generalization capability of low-thrust problems.}
\centering
\begin{tabular}{lccccc}
\hline
&\emph{m}$_s$&|$\Delta$\emph{a}|&|$\Delta$\emph{e}|&|$\Delta$\emph{i}|&|$\Delta$\emph{L}|\\\hline
Mean&2188.5 kg&0.3279 AU&0.1604&3.28 deg&32.96 deg\\
TOP&98.63\%&99.29\%&97.37\%&99.49\%&99.41\%\\
FOP&97.79\%&99.33\%&97.42\%&99.58\%&99.49\%\\
\hline
\end{tabular}
\end{table}

\begin{table}[hbt!]
\caption{\label{tab_gc_imp} Generalization capability of multi-impulse problems.}
\centering
\begin{tabular}{lccccc}
\hline
&|$\Delta$\emph{a}|&|$\Delta$\emph{e}|&|$\Delta$\emph{i}|&|$\Delta\Omega$|\\\hline
Mean&300.68 km&0.0209&3.95 deg&4.17 deg\\
JMP&94.25\%&92.10\%&96.11\%&96.54\%\\
\hline
\end{tabular}
\end{table}

\subsection{Applications in Multitarget Mission Design}

Here, the well-trained DNNs are applied in the design of three multitarget missions.
The first two multitarget missions are  multitarget main-belt asteroid, low-thrust rendezvous missions. The asteroids' orbital parameters are from the 7th GTOC \cite{casalino2014problem}. 
The third multitarget mission is a multitarget ADR J$_2$-perturbed multi-impulse mission. The debris orbital parameters are from the 9th GTOC \cite{izzo2018kessler}.

In the first mission, 10 asteroids must be visited in the minimum time. The modified Julian date (MJD) of departure is 59,215. After excluding asteroids with eccentricities larger than 0.1, 4,986 asteroids remain. The amount of all possible solutions is 4,986!/4,976! (approximately 10$^{37}$). Transfers between the two targets are time-optimal, low-thrust transfers, and the optimal transfer time and mass consumption are estimated by the well-trained DNN. Beam search, which is commonly used for this type of problem 
\cite{li2017j2,izzo2016designing}, is applied to search the optimal sequence. The same pre-screening procedure is adopted as described in Algorithm \ref{alg2}. The beam width is set to be 10,000, and the heuristic cost is set to be transfer time. One can tune the options of beam search, but it is not necessary in the present case.

After the searching process is over, the best 20 sequences are optimized using the indirect method. Results are shown in Table \ref{tab_top}, where the actual solutions are solved by the indirect method and estimated solutions are acquired by the DNN. Most of the accuracies are higher than 97\%, which is quite a high accuracy considering that the length of the sequence is 10.	The forward sequences in the actual rank are still forward in the estimated rank. The rank difference is due to the fact that their transfer days are very close to each other.

In the second mission, 10 asteroids must be visited using minimum fuel consumption. Transfers between the two targets are fuel-optimal, low-thrust transfers. The minimum transfer time $t_{min}$ is first estimated by the DNN and the transfer time is given by $t_{fop}=c*t_{min}$. In this simulation, c is set to be 1.35. The fuel consumption is then estimated by the DNN. Beam search is also applied. Other simulation parameters are the same as in the first mission. After the searching process is over, the best 20 sequences are optimized using the indirect method. Results are shown in Table \ref{tab_fop}. All the accuracies are higher than 96\%, and the actual and estimated ranks have a high consistency.

In the third mission, 100 pieces among the total 123 debris of debris must be removed  with 10 spacecraft. Each spacecraft will remove 10 debris using minimum $\Delta$\emph{v}$_{total}$. For each spacecraft, the start and transfer times are randomly given, and a pre-screening procedure, described in Algorithm \ref{alg3}, is also utilized. The single leg $\Delta$\emph{v} is estimated by the DNN. After 10 sequences are obtained, PSO is applied to optimize the actual $\Delta$\emph{v}$_{total}$. Results are shown in Table \ref{tab_imp}. All the accuracies are higher than 96\%. It is interesting that even if the accuracy of the well-trained DNN of JMPs is not as good as that of TOPs and FOPs, an excellent performance in the sequence accuracy is still evident. Results of the three multitarget missions demonstrate that DNNs can provide highly accurate solutions in preliminary mission design.

Another important factor in multitarget mission design is computation time, which is shown in Table \ref{tab_time}.  The computation time is single-leg computation time calculated by averaging the computation times when optimizing the entire sequence. The actual computation time includes the procedure that solves the problem 10 (for TOPs and FOPs) or three (for JMPs) times and chooses the best solution. There is no doubt that DNNs have an obvious advantage in computation speed.

\clearpage

\begin{table}[hbt!]
\caption{\label{tab_top} Results of time-optimal multitarget mission.}
\centering
\begin{tabular}{ccccccccccccccc}
\hline
Actual rank&Rendezvous sequence&$t_{estimated}$ (d)&$t_{actual}$ (d)&Accuracy (\%)&Estimated rank\\\hline
1&\tabincell{c}{	6553	7400	803	12570	5701\\11047	2053	3427	11005	693	}&	1287.49& 	1321.36& 	97.44& 	1\\
2&\tabincell{c}{	6553	7400	803	12570	5701\\	11047	2053	3427	693	11005}&	1297.48& 	1322.16& 	98.13& 	9\\
3&\tabincell{c}{	7424	803	7400	12570	5701\\	11047	2053	3427	11005	693}&	1288.25& 	1323.62& 	97.33& 	2\\
4&\tabincell{c}{	7424	803	7400	12570	5701\\	11047	2053	3427	693	11005}&	1298.26& 	1325.12& 	97.97& 	11\\
5&\tabincell{c}{	7424	803	7400	11047	12570\\	5701	2053	3427	11005	693}&	1293.04& 	1330.66& 	97.17& 	6\\
6&\tabincell{c}{	6553	803	7400	12570	5701\\	11047	2053	3427	11005	693}&	1292.16& 	1330.75& 	97.10& 	4\\
7&\tabincell{c}{	6553	7400	803	11047	12570\\	5701	2053	3427	11005	693}&	1290.68& 	1331.59& 	96.93& 	3\\
8&\tabincell{c}{	6553	803	7400	11047	12570\\	5701	2053	3427	11005	693}&	1295.14& 	1332.15& 	97.22& 	8\\
9&\tabincell{c}{	7424	803	7400	11047	12570\\	5701	2053	3427	693	11005}&	1303.04& 	1334.50& 	97.64& 	17\\
10&\tabincell{c}{	6553	803	7400	12570	5701\\	11047	2053	3427	693	11005}&	1302.17& 	1334.61& 	97.57& 	14\\
11&\tabincell{c}{	7424	803	7400	11047	5701\\	12570	2053	3427	11005	693}&	1302.96& 	1335.33& 	97.58& 	16\\
12&\tabincell{c}{	6553	7400	803	11047	12570\\5701	2053	3427	693	11005}	&	1300.70& 	1335.70& 	97.38& 	13\\
13&\tabincell{c}{	6553	803	7400	11047	12570\\	5701	2053	3427	693	11005}&	1305.12& 	1336.42& 	97.66& 	20\\
14&\tabincell{c}{	6553	7400	803	11047	5701\\12570	2053	3427	11005	693}	&	1298.34& 	1336.58& 	97.14& 	12\\
15&\tabincell{c}{	803	7400	6553	12570	5701\\	11047	2053	3427	11005	693}&	1293.17& 	1337.27& 	96.70& 	7\\
16&\tabincell{c}{	803	7400	6553	12570	5701\\	11047	2053	3427	693	11005}&	1303.17& 	1342.70& 	97.06& 	18\\
17&\tabincell{c}{	7400	803	6553	12570	5701\\	11047	2053	3427	11005	693}&	1292.47& 	1368.52& 	94.44& 	5\\
18&\tabincell{c}{	3914	9641	10952	4932	3013\\8433	7359	8484	7366	9616}&	1297.95& 	1370.37& 	94.72& 	10\\
19&\tabincell{c}{	7400	803	6553	12570	5701\\	11047	2053	3427	693	11005}&	1302.48& 	1371.74& 	94.95& 	15\\
20&\tabincell{c}{	6553	7400	803	12570	11047\\	5701	2053	3427	11005	693}&	1304.24& 	1378.53& 	94.61& 	19\\
\hline
\end{tabular}
\end{table}

\clearpage

\begin{table}[hbt!]
\caption{\label{tab_fop} Results of fuel-optimal multitarget mission.}
\centering
\begin{tabular}{ccccccccccccccc}
\hline
Actual rank&Rendezvous sequence&$m_{estimated}$ (kg)&$m_{actual}$ (kg)&Accuracy (\%)&Estimated rank\\\hline
1&\tabincell{c}{1150	14205	7410	12512	12461\\	4984	14195	5033	3018	3911 }	&	1412.47& 	1376.76& 	97.41&1\\ 
2&\tabincell{c}{	1150	14205	12512	7410	12461\\	4984	14195	5033	3018	3911}&	1401.73& 	1362.17& 	97.10&2 \\ 
3&\tabincell{c}{	14205	1150	7410	12512	12461\\	4984	14195	5033	3018	717}&	1382.29& 	1360.13& 	98.37&8 \\ 
4&\tabincell{c}{	1150	14205	7410	12512	12461\\	4984	14195	3018	5033	717}&	1382.90& 	1359.52& 	98.28&7 \\ 
5&\tabincell{c}{	1150	14205	7410	12512	4984\\	12461	14195	5033	3018	717}&	1393.23& 	1359.36& 	97.51&4 \\ 
6&\tabincell{c}{	1150	14205	7410	12512	4984\\	12461	14195	5033	3018	3911}&	1393.51& 	1354.31& 	97.11&3 \\ 
7&\tabincell{c}{	14205	1150	7410	12512	12461\\	4984	14195	5033	3018	3911}&	1382.93& 	1354.13& 	97.87&6 \\ 
8&\tabincell{c}{	1150	14205	7410	12512	12461\\	4984	14195	5033	717	7362}&	1381.12& 	1351.52& 	97.81&9 \\ 
9&\tabincell{c}{	1150	14205	7410	12512	12461\\	4984	7393	7362	717	5033}&	1377.29& 	1350.86& 	98.04&13 \\ 
10&\tabincell{c}{	1150	14205	7410	12512	12461\\	4984	7362	7393	3911	5033}&	1374.63& 	1350.83& 	98.24&16 \\ 
11&\tabincell{c}{	14205	1150	12512	7410	12461\\	4984	14195	5033	3018	717}&	1374.08& 	1349.14& 	98.15& 18\\ 
12&\tabincell{c}{	6534	14205	12512	7410	12461\\	4984	14195	5033	3018	717}&	1383.67& 	1348.34& 	97.38&5 \\ 
13&\tabincell{c}{	1150	14205	7410	12512	12461\\	4984	14195	5033	3018	7362}&	1377.42& 	1347.33& 	97.77&12 \\ 
14&\tabincell{c}{	1150	14205	7410	12512	4984\\	12461	14195	5033	717	3018}&	1378.52& 	1343.63& 	97.40&10 \\ 
15&\tabincell{c}{	1150	14205	7410	12512	12461\\	4984	14195	3018	3911	5033}&	1375.53& 	1342.19& 	97.52&14 \\ 
16&\tabincell{c}{	1150	14205	7410	12512	12461\\	4984	7362	7393	717	5033}&	1373.63& 	1340.00& 	97.49&19 \\ 
17&\tabincell{c}{	6534	14205	7410	12512	4984\\	12461	14195	5033	3018	717}&	1373.53& 	1337.48& 	97.30&20 \\ 
18&\tabincell{c}{	1150	14205	7410	12512	12461\\	4984	14195	5033	3018	2645}&	1377.83& 	1337.13& 	96.96&11 \\ 
19&\tabincell{c}{	1150	14205	12512	7410	12461\\	4984	14195	3018	5033	3911}&	1374.71& 	1336.06& 	97.11&15 \\ 
20&\tabincell{c}{	6534	14205	7410	12512	4984\\	12461	14195	5033	3018	3911}&	1374.27& 	1332.29&	96.85&17 \\ 

\hline
\end{tabular}
\end{table}

\begin{table}[hbt!]
\caption{\label{tab_imp} FOP.}
\centering
\begin{tabular}{ccccccccccccccc}
\hline
Spacecraft&Removal sequence&$\Delta\emph{v}_{DNN}$ (m/s)&$\Delta\emph{v}_{actual}$ (m/s)&Accuracy (\%)\\\hline
1&31 1 5 78 24 23 38 11 32 16&	768.00& 	798.80& 	96.14\\ 
2&52 8 21 18 62 41 9 108 115 59&	965.29& 	965.39& 	98.86\\
3&4 63 37 26 109 79 20 61 15 112&	1288.36& 	1285.29& 	99.76\\
4&100 29 120 86 55 14 67 46 53 94&	1038.73& 	1014.49& 	97.61\\
5&106 85 93 73 39 36 45 50 95 98&	1072.01& 	1064.22& 	99.27\\
6&10 91 99 65 22 70 71 33 7 13&	1436.10& 	1434.62& 	99.90\\
7&68 102 77 34 54 118 3 113 47 81&	1194.01& 	1195.71& 	99.86\\
8&111 58 114 88 27 76 110 101 80 105&	719.16& 	730.01& 	98.51\\
9&119 42 0 60 40 43 12 66 83 6&	996.17& 	1012.42& 98.40	\\
10&30 17 57 92 56 69 97 2 72 87&	908.91& 	925.45& 98.21	\\
\hline
\end{tabular}
\end{table}

\begin{table}[hbt!]
\caption{\label{tab_time} Computation time.}
\centering
\begin{tabular}{ccccccccccccccc}
\hline
&time$_{DNN}$ (ms)&time$_{actual}$ (ms)\\\hline
TOP&6.1&2061&\\ 
FOP&3.5&17337&\\ 
JMP&3.2&399782&\\ 
\hline
\end{tabular}
\end{table}

\section{Conclusions}

In this paper, three types of optimal transfers, time-optimal low-thrust transfers, fuel-optimal low-thrust transfers, and minimum-$\Delta$\emph{v} J$_2$-perturbed multi-impulse transfers, are estimated using DNNs. Random search is recommended to search hyper-parameters of DNNs because it is highly effective and easy to implement. After the DNNs are well-trained, the accuracy can be higher than 99.5\% for low-thrust trajectories and higher than 96\% for multi-impulse trajectories. DNNs are also compared with other machine-learning algorithms that have already been proved to improve accuracy in the estimation of optimal transfers, and are found to be superior to all of them. The investigations of generalization capability show that well-trained DNNs can still estimate optimal transfers accurately even if transfer parameters are outside the area in which the datasets are generated. Finally, DNNs are applied in three multitarget mission designs and are proved to be highly fast and accurate. This work demonstrates the estimation of optimal transfers via a DNN trained on a regression task, and further applications in multitarget mission design can be sought.

\section*{Funding Sources}
This research is supported by the Chinese National Natural Science Fund for Distinguished Young Scientists of China (No. 11525208).

\bibliography{ref}

\end{document}